\documentclass{article} 
\newcommand{\non}{\nonumber}
\newcommand{\s}{\\[1ex]}
\newcommand{\beq}{\begin{equation}} 
\newcommand{\eeq}{\end{equation}}
\newcommand{\lbl}{\label} 
\newcommand{\q}{\quad} 
\newcommand{\qq}{\quad\quad} 
\newcommand{\tr}{{\rm tr\,}} 
\newcommand{\arrow}{\rightarrow} 
\newcommand{\cs}{{\cal S}} 
\newcommand{\goesto}{{\scriptstyle \rightarrow}}
\newcommand{\goesupto}{{\scriptstyle \uparrow}} 
\newcommand{\re}[1]{(\ref{#1})}
\newcommand{\rang}{\mbox{\rm range\,}} 
\newcommand{\rnk}{\mbox{\rm rank\,}} 
\newcommand{\ind}{\mbox{index\,}} 
\newcommand{\kr}{\mbox{ker\,}} %Note:  control sequence \ker is apparently
\newcommand{\coker}{\mbox{coker\,}} 
\newcommand{\qedd}{ \hfill \vrule height4pt width3pt depth2pt \vskip .5cm}
\newenvironment{proof}{\bf Proof: \rm}{\qedd}
\newenvironment{namedproof}[1]{\noindent \bf #1 \rm}{\qedd}

\newtheorem{theorem}{Theorem}
\newtheorem{lemma}[theorem]{Lemma}
\newtheorem{corollary}[theorem]{Corollary}
\newtheorem{proposition}[theorem]{Proposition}
 
\newcommand{\zbar}{\bar{z}} 
\newcommand{\cl}{{\cal L}} 
\newcommand{\cls}{{{\cal L}_*}} 
\newcommand{\clp}{{{\cal L}^\prime}} 
\newcommand{\ml}{{M(\cl)}}
\newcommand{\mlp}{{M(\clp)}} 
\begin{document}
\title{The Curvature of a Single Contraction Operator on a Hilbert Space 
\protect{
\thanks{AMS Subject Classification: 47A13 (Primary); 47A20 (Secondary).
Keywords:  operator, curvature }
} 
}
\author{Stephen Parrott \\ 
Department of Mathematics and Computer Science \\
University of Massachusetts at Boston \\
100 Morrissey Blvd.\\
Boston, MA 02125\\
USA  
\protect{
\thanks{This work was done while on sabbatical leave at the University
of California at Berkeley.  I would like to express my gratitude 
for the hospitality of its Mathematics Department, both current and in
years past.} 
}
\date{June 29, 2000}
}
\maketitle
\begin{abstract}
This note studies  
Arveson's curvature invariant for $d$-contractions
$T = (T_1, T_2, \ldots , T_d)$ for the special case $d = 1$, 
referring to a single contraction operator $T$ on a Hilbert space.  
It establishes a formula which gives an easy-to-understand
meaning for the curvature of a single contraction. 
The formula is applied to give an example of an operator with nonintegral
curvature. 
Under the additional hypothesis that the single contraction $T$ 
be ``pure'', we show that its curvature $K(T)$ 
is given by  $ K(T) = - \ind (T) := - (\dim \kr (T) - \dim \coker (T))$. 
\end{abstract} 
\noindent 
\section{The curvature of a single operator}
This note studies  
Arveson's curvature invariant for $d$-contractions
$T = (T_1, T_2, \ldots , T_d)$ for the special case $d = 1$, 
referring to a single contraction operator $T$ on a Hilbert space.  
It establishes a formula which gives an easy-to-understand
meaning for the curvature of a single contraction. 
The formula is applied to give an example of an operator with nonintegral
curvature. 
Under the additional hypothesis that the single contraction $T$ 
be ``pure'', 
we show that its curvature $K(T)$ (defined below) 
is given by  $ K(T) = - \ind (T) := - (\dim \kr (T) - \dim \coker (T))$. 

Let $T$ be a contraction operator on a Hilbert space $H$,
and $ \Delta_T := \sqrt{1-TT^*}$.
Assume that $\Delta_T$ has finite rank.
Then the curvature $K(T)$ of $T$ (our shorthand for what should
properly be called
the  curvature of the Hilbert module associated with $T$)
is defined in \cite{arvproc} as:
\beq
\lbl{curvdef}
K(T) := \int_{|z| = 1}  dz \, \lim_{r \goesupto 1} \, (1-r^2) 
\, \tr (\Delta_T (1-rzT^*)^{-1} (1-r\zbar T)^{-1} \Delta_T)
\q.
\eeq

This is a specialization to the case of a single operator of Arveson's more 
general theory of $d$-contractions, 
which are finite sets of $d$ commuting operators 
satisfying an auxiliary condition analogous to contractiveness of our $T$. 

We refer the reader to \cite{arvsub}, \cite{arvproc},  and \cite{arvprep}
for the definition and basic properties of $d$-contractions.
However, we consider $d$-contractions for $d > 1$ solely 
for purposes of placing our results 
within the framework of the more general theory, 
and essentially no knowledge of $d$-contractions is necessary 
to follow our proofs.
The only reliance on the general theory is that
Arveson's Stability of Curvature result, \cite{arvprep},
Section 3, Corollary 1, is used in the proof  
of Proposition \ref{prop1}.  However, as noted there,
the reader can easily establish 
this result directly for the special case of a 1-contraction,
which is all that we need. 

The definition of curvature implicitly assumes the existence of the
limit in \re{curvdef}.  A theorem stated in \cite{arvproc},
and proved in \cite{arvprep} (Theorem A), 
guarantees the existence of the limit for almost all $z$,
and moreover bounds it above by the rank of $\Delta_T$.
For the case of a single operator, 
this also follows from the discussion  
of \cite{nagy-foias}, Chapter VI, Section 1, particularly,
page 238, equation (1.5).   

Let $ T: H \arrow H$ be a contraction on a Hilbert space $H$ with 
$\rnk \sqrt{1 - TT^*} $ finite.  Note that this implies that
$\rang \sqrt{1 - TT^*} = \rang (1 - TT^*) $, 
a fact which will be used frequently without comment.

First we associate with $T$ a partial isometry $Q$
with the same curvature, so that for most purposes of computing the curvature,
we may assume that $T$ is itself a partial isometry.
This is not always necessary, 
but it makes many problems easier to think about. 
\begin{proposition} 
\label{prop1}
With $T$ as just described, set
$$
Q := 
\left[
\begin{array}{cc}
T & \sqrt{1-TT^*} \\
0 & 0
\end{array}
\right]
\q,
$$
considered as an operator on $H \oplus \rang (1-TT^*)$. 
\s
Then $Q$ is a partial isometry with $K(Q) = K(T)$.
\s
Moreover, $\rnk (1-QQ^*) = \rnk(1-TT^*)$,
\\
and,
 $\rnk (1-Q^*Q) = \rnk(1-T^*T)$. 
\end{proposition}
\begin{proof} 
That  $K(Q) = K(T)$ follows from one of Arveson's key results 
for $d$-contractions,
Stability of Curvature, \cite{arvprep}, Section 3, Corollary 1. 
For our case of a 1-contraction, a proof can alternatively  
be obtained by a straightforward calculation of $K(Q)$,
based on its definition \re{curvdef}.  

Since
\beq
\lbl{Qcoker}
1 - QQ^* = 
\left[
\begin{array}{ll}
0 & 0 \\ 
0 & 1_{\mbox{\scriptsize Range\,} (1-TT^*)} 
\end{array}
\right]
\q,
\eeq
it is obvious that $Q$ is a partial
isometry with $\rnk (1-QQ^*) = \rnk (1 - TT^*)$.  

Next we show that $\rnk (1-Q^*Q) = \rnk (1-T^*T)$. 
We have
$$ 
1 - Q^*Q = 
\left[
\begin{array}{cc}
1 - T^*T & -T^*\sqrt{1-TT^*} \\
-\sqrt{1-TT^*}T & TT^* 
\end{array}
\right]
\q.
$$
Since the rank of an operator matrix is at least as large as the rank
of any entry, 
if $\rnk (1-T^*T)$ is infinite, so is  $\rnk (1 - Q^*Q)$.
Thus we may assume that $\rnk(1-T^*T)$ is finite.

Let $C_i,\  \, i = 1,2$,  denote the $i$'th column of the matrix 
for $1-Q^*Q$, considered in the obvious way as operators, e.g.,
$C_1 : H \arrow \rang  (1-TT^*)$.
Then
$$
C_2 \sqrt{1-TT^*} = 
\left[ 
\begin{array}{c}
-T^*(1-TT^*) \\
TT^*\sqrt{1-TT^*}
\end{array}
\right] 
=  
\left[ 
\begin{array}{c}
-(1-T^*T)T^* \\
\sqrt{1-TT^*}TT^*
\end{array}
\right] 
= - C_1 T^*
$$
Since the domain of $C_2$ is $\rang \sqrt{1-TT^*}$, 
this implies that $\rang C_2 \subset \rang C_1$, and hence
$\rang (1-Q^*Q) = \rang C_1$. 

It is well known (e.g., \cite{riesz-nagy}, Section 147) that 
$\sqrt{1-TT^*}\, T = T \sqrt{1-T^*T}$, so
$$
C_1 =  
\left[ 
\begin{array}{c}
(1-T^*T) \\
-T\sqrt{1-T^*T} 
\end{array}
\right] 
= 
\left[ 
\begin{array}{c}
\sqrt{1-T^*T} \\
-T
\end{array}
\right] 
\sqrt{1-T^*T}
$$
Since 
$$
\left[ 
\begin{array}{c}
\sqrt{1-T^*T} \\
-T
\end{array}
\right] 
$$
is an isometry, the map $\sqrt{1-T^*T} x \mapsto C_1 x  , \ x \in H$,
defines an isometric bijection between 
$\rang \sqrt{1-T^*T}$ and $\rang \, C_1 = \rang (1-Q^*Q)$.  
Hence $\rnk (1-T^*T) = \rnk (1-Q^*Q)$.
\end{proof}

Next we derive a simple formula for $K(Q)$, along with a variant formula
for $K(T)$ which does not mention $Q$.   
The formula for $K(Q)$ seems particularly helpful in thinking
about these problems.
\begin{theorem}
\label{curvform}
Let $Q$ be a partial isometry such that $\Delta_Q := \sqrt{1-QQ^*}$
has finite rank, and 
let $e_1, e_2, \ldots , e_q$ be an orthonormal basis for 
$\rang(\Delta_Q) $. Then
$$
K(Q) = \sum_{k=1}^q \lim_{n \goesto \infty} \| Q^n e_k \|^2
\q.
$$
Moreover, for any contraction $T$ 
for which $\Delta_T$ has finite rank,
$$
K(T) = \lim_{n \goesto \infty} \tr ({T^*}^n T^n (1 - TT^*)) 
\q.
$$
\end{theorem}
\begin{proof}

The boundedness of the integrand of the curvature 
justifies application of the Lebesgue Dominated Convergence Theorem 
to interchange limit and integral in the definition \re{curvdef}  
of curvature:
\begin{eqnarray}
K(Q) &:=& 
\lbl{line1}
\int_{|z|=1} \lim_{r \goesupto 1}\,  (1-r^2) 
\sum_{k=1}^q 
\langle (1-r zQ^*)^{-1} e_k , (1-r \zbar Q)^{-1} \rangle e_k \rangle 
\, dz \\
&=&
\lbl{line2}
\sum_{k=1}^q 
\lim_{r \goesupto 1} \, (1-r^2) 
\int_{|z|=1}  
\langle (1- r \zbar Q)^{-1} e_k , (1- r \zbar Q)^{-1} e_k \rangle \, dz \\
&=& 
\lbl{line3}
\sum_{k=1}^q 
\lim_{r \goesupto 1} \, (1-r^2) 
\int_{|z|=1} 
\langle \, \sum_{i = 0}^\infty (r \zbar Q)^i e_k, 
~ \sum_{j = 0}^\infty (r \zbar Q)^j e_k  \, \rangle \, dz \\ 
&=& 
\lbl{line5}
\sum_{k=1}^q 
\lim_{r \goesupto 1} \, (1-r^2) 
\int_{|z|=1} 
\sum_{i, j = 0}^\infty r^{i+j} z^{j-i} 
\langle Q^i e_k,~ Q^j e_k \rangle \, dz \\
&=& 
\lbl{line6}
\sum_{k=1}^q 
\lim_{r \goesupto 1} \, (1-r^2) 
\sum_{i = 0}^\infty r^{2i} \,  
\| Q^i e_k \|^2 \\
&=&
\lbl{line7}
\sum_{k=1}^q 
\lim_{i \goesto \infty}
\| Q^i e_k \|^2 
\end{eqnarray}

Equation \re{line6} was obtained by interchanging the infinite sum and integration.
This is justified because for fixed $r$, the infinite sum 
converges absolutely with sum of absolute values bounded above 
by $(1-r^2)^{-2}$. 

Equation \re{line7} is justified as follows. 
For fixed $k$, consider the decreasing sequence 
$$ 1 \geq \| Qe_k \| \geq \| Q^2 e_k \| \geq  
 \ldots \geq \lim_{i \goesto \infty} \|Q^i e_k \|  \q, $$  
and set $L := \lim_{i \goesto \infty} \|Q^i e_k \|^2$.
Then for any positive integer $m$,  
\begin{eqnarray*}
L &=&  L \, \lim_{r \goesupto 1} \, (1-r^2) \sum_{i=0}^\infty r^{2i} \\
&\leq& \lim_{r \goesupto 1} \, (1-r^2) 
    \sum_{i=0}^\infty r^{2i} \, \| Q^i e_k\|^2 \\
&=& \lim_{r \goesupto 1} \, (1-r^2) \sum_{i = m}^\infty r^{2i} 
\| Q^i e_k \|^2 \\
&\leq & \|Q^m e_k \|^2  \, \lim_{r \goesupto 1} \, (1-r^2) 
\sum_{i = m}^\infty r^{2i}  \\
&=& \| Q^m e_k \|^2 
\q.
\end{eqnarray*} 
For sufficiently large $m$, the right side is arbitrarily close to $L$,
showing that
$$ 
\lim_{r \goesupto 1} \, (1-r^2) \sum_{i = 0}^\infty r^{2i} 
\| Q^i e_k \|^2 = 
\lim_{i \arrow \infty} \| Q^i e_k \|^2  
\q,
$$ 
thus proving \re{line7}.

This proves the asserted formula for $K(Q)$.  To prove 
the alternative formula for $K(T)$, define $Q$ to be the 
partial isometry of Proposition \ref{prop1}
with $K(Q) = K(T)$.   Recall from \re{Qcoker} that 
$$
\Delta_Q = 
\left[
\begin{array}{cc}
0 & 0 \\
0 & 1 
\end{array}
\right]
\q,
$$
and check that 
\beq 
Q^n = 
\left[
\begin{array}{cc}
T^n & T^{n-1}\Delta_T \\
0 & 0
\end{array}
\right] 
\q.
\eeq
The formula for $ K(T)$ follows immediately upon combining these observations, 
Proposition  \ref{prop1}, 
the formula just proved for $K(Q)$, and the cyclic
property of the trace:
\begin{eqnarray*} 
K(T) &=& K(Q) =  
 \sum_{k=1}^q \lim_{n \goesto \infty} \| Q^n e_k \|^2 \nonumber \\
&=&
 \sum_{k=1}^q \lim_{n \goesto \infty} 
\langle \Delta_T {T^*}^n T^n \Delta_T e_k , e_k \rangle \nonumber \\
&=& \lim_{n \goesto \infty} \tr (\Delta_T {T^*}^n T^n \Delta_T ) \nonumber \\
&=& \lim_{n \goesto \infty} \tr ({T^*}^n T^n (1 - TT^* )) 
\end{eqnarray*}
\end{proof} 
A simple sufficient condition for the curvature to vanish
is an immediate corollary:
\begin{corollary}
Any contraction $T$ whose positive powers  $T^n$ converge strongly to 0
has vanishing curvature:  $K(T) = 0$. 
\end{corollary}
\section{Relation to Arveson's curvature formula}
Arveson \cite{arvproc} established 
a different formula for the curvature
of a $d$-contraction $T$.  Specialized to the case $d=1$, it reads:  
\beq
\lbl{arvformula}
K(T) = \lim_{n \goesto \infty} 
\frac{\tr( 1 - T^n {T^*}^n)}{n}
\q.
\eeq 

In order to make clear how our formula fits into Arveson's framework, 
we now derive ours assuming his. 
However, the resulting proof is not notably simpler than the direct
proof above,
and Arveson's proof is even more involved, corresponding to
the fact that the case $d > 1$ is probably 
fundamentally more difficult than $d = 1$.
For the single operator case $d = 1$,
Arveson's formula follows similarly from ours.

Let $T$ be a contraction with $\Delta_T$ of finite
rank, and $e_1, \ldots , e_q$ an orthonormal basis for 
$\rang(\Delta_T) = \rang (1 - TT^*)$. 
First note the collapsing sum: 
$$
1 - T^n {T^*}^n = \sum_{i = 0}^{n-1} T^i (1 - TT^*) {T^*}^i \\ 
\q.
$$
Hence
\begin{eqnarray*}
\lim_{n \goesto \infty} \frac{\tr (1 - T^n {T^*}^n)}{n} &=&  
\lim_{n \goesto \infty} 
\frac{1}{n} \sum_{i = 0}^{n-1} \tr( {T^*}^i T^i (1 - TT^*) ) \nonumber\\
&=& \lim_{i \goesto \infty}  \tr ( {T^*}^i T^i (1 - TT^*))  \q.\\ 
\end{eqnarray*} 

The last equality was obtained as follows.
Consider the sequence 
\begin{eqnarray*}
a_i &:=& \tr ( {T^*}^i T^i (1 - TT^*)) \nonumber \\
&=& \tr ( (T^i \Delta_T)^* (T^i \Delta_T)) \nonumber \\
&=& \sum_{k=1}^q \langle\  (T^i \Delta_T)^* (T^i \Delta_T)) e_k, e_k 
\rangle \nonumber \\
&=& \sum_{k=1}^q \| T^i \Delta_T e_k \|^2  
\q.
\end{eqnarray*}
The last expression makes clear that 
$ a_1 \geq a_2 \geq \ldots \geq  0$,
so that the sequence has a limit $L = \lim_{i \goesto \infty} a_i$. 
We shall show that for any such sequence $a_i$,  
$$\lim_{n\goesto\infty} (1/n) \sum_{i=0}^{n-1} a_i = L \q.$$ 
Since the sequence  
$\left\{ \frac{1}{n} \sum_{i=0}^{n-1} a_i \right\}$ is bounded above
by $a_0$, it is enough to show that its only possible accumulation 
point is $L$.

For any fixed $m$ and all $n \geq m$, 
$$
L \leq  \frac{1}{n} \sum_{i=0}^{n-1} a_i \leq 
\frac{1}{n} \sum_{i-0}^{m-1} a_i +  \frac{n-m}{n} a_m
\q.
$$
Letting $n$ tend to infinity with $m$ fixed, 
we see that any accumulation point of the sequence 
$\left\{ \frac{1}{n} \sum_{i=0}^{n-1} a_i \right\}$ must lie 
between $L$ and $a_m$.  
Finally, letting $m$ tend to infinity 
shows that $L$ is the only accumulation point.
\q.  
\section{A simple formula for the curvature of a single, pure, contraction.} 
A single contraction $T$ on a Hilbert space $H$ 
will be called {\em pure} if for all $h \in H$, $\lim_{n \goesto \infty}
{T^*}^n h = 0$; i.e., if the adjoint powers ${T^*}^n$ converge strongly
to 0.  This is the specialization to the case $d = 1$ of Arveson's more 
complicated definition of a pure $d$-contraction.   

Arveson remarked \cite{arvproc} that it is generally difficult 
to determine the curvature of a $d$-contraction,
but that ``in the few cases where the computations 
can be explicitly carried out,
the curvature turns out to be an integer.'' 
This led him to ask \cite{arvproc} 
if the curvature of a pure $d$-contraction
need always be an integer.  

This was a surprising suggestion, because 
nothing in the definition of curvature suggests  
that it should be an integer.  
Subsequently, 
D. Greene, S. Richter, and C. Sundberg \cite{greene} proved that indeed
the curvature of any pure $d$-contraction is an integer. 
However, their function-theoretic methods do not seem to give
an effective procedure for calculating this integer in particular cases, 
and a geometric, operator-theoretic interpretation
of the curvature of a general $d$-contraction remains elusive
as of this writing.  

Our contribution toward understanding the meaning of the curvature invariant
is a simple, usually easily computable, formula 
for the curvature of single, pure contraction; i.e., the special
case $d = 1$. It states that the curvature is the difference 
of the dimensions of two subspaces, and hence is obviously
integral.
The methods of proof are operator-theoretic,
based on unitary dilation theory as set forth in 
\cite{nagy-foias}.  It uses neither the 
Greene/Richter/Sundberg result nor their function-theoretic
methods, and thus gives an independent proof of their result 
for the special case $ d = 1$.  

Our characterization of the curvature of a single pure contraction is: 
\begin{theorem}
\label{curvthm} 
Let $T$ be a {\em pure} contraction operator
such that $\Delta_T := \sqrt{1 - T T^*}$ has finite rank.
Then its curvature $K(T)$ is the integer
\beq
\lbl{curvthmconc}
 K(T) = \dim \rang (1 - TT^*) - \dim \rang (1 - T^* T) 
\q.
\eeq 
\end{theorem} 
A counterexample in the next section uses Theorem \ref{curvform}
to show that the hypothesis that
$T$ be pure is essential.  
 
Before proving the theorem,
we review some standard facts about unitary dilations.
Proofs can be found
in \cite{nagy-foias}, particularly Chapters 1, 2, and 6.  
We give specific references from this work for key  facts
required by the proof.

Let $T$ be a contraction on a Hilbert space $H$, and $U$
its minimal unitary dilation to a larger Hilbert space $K \supset H$.
This means that $P_H U^n | H = T^n$ for all $n \geq 0$,
where $P_H$ denotes the projection to $H$, and minimality means 
that $K = \bigvee^\infty_{n = -\infty} U^n H $.
\begin{enumerate}
\item
\lbl{item1} 
The minimal unitary dilation $U$ for $T$ may be constructed as follows.
Define 
\beq
\lbl{kdef}
K := \ldots \oplus \overline{\Delta_{T} H} 
\oplus \overline{\Delta_{T} H} \oplus H 
\oplus \overline{\Delta_{T^*} H} \oplus \overline{\Delta_{T^*} H} \ldots ,
\eeq
where the overscore denotes closure.
(The closures turn out to be unnecessary in our context,
but that only becomes apparent later.) 
Consider $H$ as embedded in $K$ in the obvious way. 
Then $U$ is defined on $K$ by:
\begin{eqnarray}
\lbl{udil1}
\lefteqn{U( \ldots ,\  b_2,\  b_1,\  b_0, \fbox{$h$}\, , a_0,\  a_1,\  a_2, 
\ \ldots ) := } 
\nonumber \\ 
&&
(\ldots ,\  b_2,\  b_1, \fbox{$Th + \Delta_T b_0$}\, ,\ 
-T^* b_0 + \Delta_{T^*} h ,\  a_0,\  a_1,\  \ldots )
. 
\end{eqnarray}
Here $b_i \in \overline{\Delta_{T}} H $, 
$a_i \in \overline{\Delta_{T^*}} H $,  and zero'th components 
(vectors in $H$) are 
distinguished by boxes.

The realization of $U$ just given is best for some purposes,
but a change of notation will bring out more clearly the features 
which will be important to us.  
Set $\cl := (U-T)H$ and $\cls := (U^* - T^*)H$.
Informally, $\cl$ is the leftmost $\overline{\Delta_{T^*} H}$ factor in \re{kdef}.
The other $\overline{\Delta_{T^*} H}$ factors are images of the leftmost
under positive powers of $U$.
Similarly, $\cls$ is the rightmost $\overline{\Delta_T H}$ factor,
and the other $\overline{\Delta_T H}$ factors are images of it under
negative powers of $U$. 

To reflect these insights, 
instead of realizing $K$ as above, think of it as follows: 
\beq
\lbl{udil2}
K \cong \ldots \oplus U^{-2}\cls \oplus U^{-1} \cls \oplus \cls
\oplus H \oplus \cl \oplus U\cl \oplus U^2 \cl \ldots
\eeq
Here $\cong$ stands for unitary equivalence.

The conceptual advantage of \re{udil2} is that it makes clear at a glance 
much of the action of $U$ on $K$.  
Unfortunately, it is awkward for the purpose 
of defining $U$ due to logical circularity.

Embedded in $U$ are two bilateral shifts which interact 
in a complicated way.  
One shifts $\cl$, and the other shifts $\cls$.
One half of each shift is transparently visible
in \re{udil2}.  For example, $U$ obviously acts as a unilateral shift
(with multiplicity $\mbox{dim}\, \cl$)
on the invariant subspace 
\beq
\lbl{mlplusdef} 
\ml^+ := \bigoplus_{n=0}^\infty U^n \cl
\q.
\eeq
Since all iterates $U^n \cl ,\ -\infty \leq n \leq \infty$ 
are easily seen to be pairwise orthogonal, also $U$ acts
as a bilateral shift on the invariant subspace 
\beq
\lbl{mldef}
\ml := \bigoplus_{n=-\infty}^\infty U^n \cl
\q,
\eeq
but the left half of this subspace, $\ml \ominus \ml^+$,
is embedded in a non-transparent way in $K$. 

A subspace $\cs$ such that the subspaces $U^n \cs$ 
are pairwise orthogonal,  
$\ -\infty < n < \infty$, is called
a {\em wandering subspace} for $U$. 
Thus $\cl$ is a wandering subspace, and so is $\cls$.
For any wandering subspace $\cs$, we'll use the notation 
$M(\cs)$ as defined in 
\re{mldef} with $\cl$ replaced by $\cs$.
\item 
\label{item2}
The contraction $T$ is pure, i.e., ${T^*}^n \arrow 0$ strongly,
if and only if  
$$
	K = M(\cls) :=  \bigoplus^\infty_{n = -\infty} U^n \cls
\q
$$
(Chap. 2, Thm. 1.1, p. 57).
\item
\lbl{item3}
If $T$ is pure, then $\dim \cl \leq \dim \cls$; equivalently,
%$\rank \Delta_{T^*} \leq \rank \Delta_T$.  
$\rnk \Delta_{T^*} \leq \rnk \Delta_T$.  
In particular, under our hypotheses that $T$ is pure with $\Delta_T$
of finite rank, also $\Delta_{T^*}$ has finite rank,
and both $\cl $ and $\cls$ are finite dimensional. 
\s
This follows from item \ref{item2} 
above combined with \cite{nagy-foias}, Chap. 1, Prop. 2.1, p. 4.
Alternatively, it can be obtained for the case that we'll need,
$\dim \cls < \infty$, 
from the Reciprocity Lemma \ref{replemma} below
with $\cl^\prime := \cls$.  Assuming temporarily that $\dim \cl$
is known to be finite, the Reciprocity Lemma applies as follows: 
$$
\dim \cl = \tr(P_\cl ) = \tr(P_\cl P_{M(\cls)} ) = 
\tr (P_{\cls} P_{M(\cl)}) \leq \tr(P_{\cls} ) = \dim \cls
\ .
$$ 
The case of an infinite-dimensional $\cl$ can be ruled out
by applying the same reasoning with $\cl$ 
replaced by finite-dimensional subspaces of $\cl$.
\item
\lbl{item4}
When $T$ is a partial isometry,
$$
U\cls = \overline{\Delta_T H} 
\q.
$$ 
In particular, $U\cls \subset H$.
\s
This is immediate from \re{udil1} after recalling that 
a partial isometry $T$
satisfies $T^* (1-TT^*)  H  =  \{0\}$. 
\end{enumerate} 

Let  $E$ and $F$ be projections on a Hilbert space, at least one
of which has finite rank. 
Then $\tr (EF) = \tr(E^2 F) = \tr (EFE)$, so $ \tr(EF)$ is 
always non-negative, is zero if and only if $E$ and $F$ have
orthogonal ranges, and takes on its maximum value 
$\mbox{dim\,($E$)}$ or 
$\mbox{dim\, ($F$)}$ only when $E \leq F$ or $F \leq E$. 
Thus  $\tr (EF) $
serves as a measure of how nearly the ranges of $E$ and $F$ coincide. 
For lack of a standard term, call $\tr (EF) $ the {\em affinity}
between the ranges of $E$ and $F$. 

The following lemma, which we call  
the Reciprocity Lemma, may have some interest in its own right.  
It states that for wandering subspaces $\cl$ and $\clp$
for a unitary operator $U$, the affinity between $\cl$ and 
the closed span of the iterates $U^n \clp ,\  -\infty \leq n  \leq \infty$
is invariant under interchange of $\cl $ and $\clp$. 
\begin{lemma}
\lbl{replemma}
{\bf [Reciprocity Lemma]}
Let $\cl$ and $\clp$ be finite dimensional wandering subspaces
for a unitary operator $U$ on a Hilbert space $K$, and set 
$$
\ml := 
\bigoplus_{n=-\infty}^\infty U^n \cl
\qq \mbox{and} \qq 
\mlp := 
\bigoplus_{n=-\infty}^\infty U^n \clp
\q.
$$
Then, denoting by $P_S$ the projection on an arbitrary subspace $S$ of $K$, 
$$
\tr (P_\cl P_{\mlp} ) = 
\tr (P_\clp P_{\ml} ) 
\q.  
$$ 
\end{lemma}
\begin{proof} 

Note that
\beq
\lbl{lemma5eq1}
P_{M(\cl)} = \sum_{n= - \infty}^\infty P_{U^n \cl} 
= \sum_{n=-\infty}^\infty U^n P_\cl U^{-n}
\q,
\eeq
the sums converging in the strong operator topology.

Multiply \re{lemma5eq1} by $P_{\clp}$ on the left, take the trace of both sides,  
and suppose we can justify
an interchange of sum and trace, obtaining
\begin{eqnarray}
\lbl{lemma5eq2}
\tr (P_\clp P_{M_\cl}) &=& 
\tr(\sum_{n=-\infty}^\infty P_\clp U^n P_\cl U^{-n}) \non \\ 
&=& \sum_{n=-\infty}^\infty \tr (P_\clp U^n P_\cl U^{-n}) 
\q.
\end{eqnarray}
Then the following simple calculation establishes the lemma:
\begin{eqnarray} 
\lbl{lemma5eq3}
\tr ( P_{\clp} P_{M(\cl)}) &=&  
 \sum_{n=-\infty}^\infty \tr (P_\clp U^n P_\cl U^{-n}) \nonumber \\ 
&=& \sum_n \tr (P_\cl U^{-n} P_\clp U^n ) \nonumber \\
&=& \tr (P_\cl P_{M(\clp)})
\q,
\end{eqnarray}
where the last line was obtained from \re{lemma5eq2} with $\cl$
and $\clp$ interchanged and the summation index $n$ replaced by $-n$. 

The interchange of sum and trace required to justify the above calculation
is not immediate because the trace is not continuous
in the strong operator topology.  
However, the trace is well-known to be a {\em normal}
linear functional, which implies that 
for any increasing sequence of trace class positive  operators 
$$
A_1 \leq A_2 \leq \ldots \leq A_m \leq \ldots 
$$
converging in the strong operator topology to a trace class
operator $A$, we have 
\beq
\lbl{lemma5eq4}
\lim_{m \goesto \infty} \tr (A_m) = \tr (A)
\q.
\eeq
This property is a slight specialization of the definition of normality.
It follows routinely from the definition  
$\tr A := \sum_{i=1}^\infty \langle Ae_i, e_i \rangle$,
with $\{e_i \} $ an orthonormal basis. 

Noting that 
$$\tr(P_\clp U^n P_\cl U^{-n}) = 
\tr(P_\clp U^n P_\cl U^{-n}P_\clp)$$
and applying \re{lemma5eq4}
with $A_m := \sum_{n=-m}^m 
\tr(P_\clp U^n P_\cl U^{-n}P_\clp)$ 
proves \re{lemma5eq2}.  

I thank W. Arveson for suggesting the above proof to replace
the unattractive direct calculation of an earlier draft. 
\end{proof} 
\begin{namedproof}{Proof of Theorem 4:} 

Proposition \ref{prop1}
shows that we may assume that $T$ is a partial isometry. 
We are going to use Theorem \ref{curvform}
to calculate $K(T)$ by calculating 
$$
\lim_{n \arrow \infty} \| T^n e \|^2
$$
for $e \in \Delta_T H$.  
For any $h \in H$, 
$$
\| T^n h \|^2 = \| P_H U^n h \|^2
\q.
$$
Since $U^n h \in H \bigoplus M(\cl)^+ = H \bigoplus (\bigoplus_{k=0}^\infty
U^k \cl$),   
\begin{eqnarray}
\lbl{norm1}
\| T^n h \|^2 &=& \| U^n h \|^2 - 
\| P_{M(\cl)^+} U^n h \|^2 \nonumber \\
&=& \|h\|^2 - \| P_{M(\cl)^+} U^n h \|^2 
\q.
\end{eqnarray}

Write $K = M(\cl) \bigoplus R$, where (as always), the direct sum
denotes an orthogonal direct sum, so this defines the subspace $R$,
which reduces $K$ because $M(\cl)$ does. 
Then any $k \in K$ can be written 
$$
k = \sum_{i=-\infty}^\infty U^i f_i + r 
\q,
$$
with $f_i \in \cl$ and $r \in R$.
And, for any $h \in H$,
\beq
\lbl{hrep}
h = \sum_{i=-\infty}^{-1} U^i f_i + r 
\q.
\eeq

Substituting \re{hrep} in \re{norm1} gives:
\begin{eqnarray} 
\lbl{normlim}
\lim_{n \goesto \infty} \| T^n h \|^2 
&=& \| h \|^2 - \lim_{n \goesto \infty} \sum_{i=-n}^{-1} \| f_i \|^2 
\nonumber \\ 
&=& \| h \|^2 - \| P_\ml h \|^2 \nonumber \\
&=& \| h \|^2 - \langle \, P_{\ml} h,\  h \rangle
\q.
\end{eqnarray}

Choose an orthonormal basis $e_1, \ldots , e_q $ for 
$\rang \, \Delta_T H$.  
Then 
substituting \re{normlim} in Theorem \ref{curvform} 
gives:
\begin{eqnarray}
\lbl{KT}
K(T) &=& \sum_{i=1}^q \lim_{n \goesto \infty} \| T^n e_i \|^2 \nonumber\\
&=& q - \sum_{i=1}^q \langle \, P_\ml e_i ,\ e_i \rangle \nonumber\\
&=& q - \tr (P_\ml P_{\Delta_{T} H }) 
\q.
\end{eqnarray} 
Item \ref{item2} remarked that for a partial isometry $T$,
$\Delta_T H = U\cls$, and substituting this in \re{KT} gives: 
\beq
\lbl{KT2}
K(T) = q - \tr (P_\ml P_{U\cls}) 
\q.
\eeq 
By the Reciprocity Lemma \ref{replemma},
$$
\tr ( P_\ml P_{U\cls}) = \tr (P_{M(U\cls)} P_\cl)
\q.
$$
But obviously, 
$$
M(U\cls) := \bigoplus_{n=-\infty}^\infty U^n U\cls 
= \bigoplus_{k=-\infty}^\infty U^k \cls = M(\cls)
\q.
$$
Combining these facts gives the desired conclusion:
\begin{eqnarray}
K(T) &=& q -  \tr ( P_\ml) P_{U\cls}) = q - \tr (P_{M(U\cls)} P_\cl) \nonumber\\ 
&=& q - \tr(P_{M(\cls)} P_\cl) = q - \tr(P_\cl) \nonumber\\ 
&=& \dim \rang\, \Delta_T - \dim \rang \Delta_{T^*} \nonumber\\
&=& \dim \rang\, (1-TT^*) - \dim \rang (1-T^*T) \nonumber 
\q.
\end{eqnarray} 
The second line follows from item \ref{item2}'s observation that
the hypothesis that $T$ be pure is equivalent to $M(\cls) = K$.  
\end{namedproof} 

Recall that a Fredholm operator $T$ 
is one with closed range and finite-dimensional kernel and cokernel
(denoted $\kr (T)$ and $\coker (T) := \ker (T^*)$).  
The {\em index} of a Fredholm
operator $T$ is defined by
\beq
\lbl{index}
\ind (T) := 
 \dim \kr (T) - \dim \coker T
\q.
\eeq 
A fundamental theorem (e.g., \cite{douglas}, p. 128, Thm. 5.36)
states that the index is invariant under compact
perturbations: 
for any Fredholm operator $T$ and compact operator $C$, 
$T+C$ is Fredholm, and 
$\ind (T + C) = \ind (T)$.  

Formula \re{index} makes sense when $T$ has finite-dimensional 
kernel and cokernel even if $T$ doesn't have closed range.
However,
since the closed range hypothesis is needed to prove the fundamental
theorem just mentioned, the term ``index'' is generally 
restricted to Fredholm operators. 
Nevertheless, for purposes of the present
exposition, it will be convenient to broaden the definition
of $\ind (T)$ to include cases in which $T$ has finite-dimensional
kernel and cokernel, but not necessarily closed range. 

When told of the curvature formula \re{curvthmconc}
given by Theorem \ref{curvthm},
W. Arveson remarked that it looked something like an operator index and 
that he had  been working on a conjecture
that under appropriate hypotheses, 
the curvature of a $d$-contraction 
would be the index of an associated operator which he calls $D_+$, 
reminiscent of the Dirac operator.  
Shortly thereafter, he wrote up these results in \cite{arvprep2},
which proves this 
for $d$-contractions whose associated Hilbert modules are 
finite rank, pure, and graded, in the terminology of \cite{arvprep}. 
It asks if the ``graded'' hypothesis can be removed, 
and also if the associated operator $D_+$ necessarily has closed range
(and so is Fredholm). 

For the case of a 1-contraction, the associated operator $D_+$
is unitarily equivalent to $T$.   
Corollary \ref{curvcor} below observes 
that under the hypotheses of Theorem \ref{curvthm},
$T$ is Fredholm,
and its curvature equals  $ - \ind (T)$.  
The differences between Corollary \ref{curvcor} 
and the specialization of Arveson's 
result to the single operator case are that the closed range property  
is proved for $d = 1$, and the ``graded''
hypothesis is not needed. 
This holds out hope that the ``graded'' and ``closed range''
hypotheses might be removable for $d$-contractions with $d > 1$.  

The interest in identifying the curvature with an index,
apart from its evident aesthetic appeal, is that the index is 
stable under compact perturbations, but the curvature is not
known to possess such stability.   
The strongest result along these lines known as of this writing
is \cite{arvprep}, Corollary 1, Stability of Curvature, 
which proves stability of the curvature 
under certain special finite rank perturbations.  
Arveson \cite{arvprep2} notes that
removing the ``closed range'' hypothesis
would establish a much stronger stability of curvature result,
and removing the ``graded'' hypothesis would strengthen it
further.   

\begin{corollary}
\lbl{curvcor}
Let $T$ be an operator satisfying the hypotheses of Theorem 4.
Then  $T$ is Fredholm, and 
$$
K(T) = - {\rm index}\, (T)
\q.
$$ 
\end{corollary} 
\begin{proof} 
Let $T$ be an operator on a Hilbert space $H$ 
satisfying the hypotheses of Theorem \ref{curvthm}. 
First we sketch the simple proof that  $T$ must be Fredholm.

The assumed finiteness of the rank of $ 1-TT^*$ implies
that $\coker (T)$ is finite-dimensional. 
Since we have already noted that $\rnk (1-T^*T) \leq \rnk (1-TT^*)$,
also $\kr (T)$ is finite-dimensional.
That $T$ must have closed range under these circumstances can be easily seen
by noting that closed range 
is equivalent to a gap above 0 in the spectrum of $T^*T$.  
If there were not such a gap, then $1-T^*T$ would not have
finite rank.

To show that \re{curvthmconc} equals $ -\ind (T)$,
let $\tilde{T}$ be the operator on $H \oplus \rang (1 - TT^*) $ defined by
the operator matrix:
$$
\tilde{T} := 
\left[
\begin{array}{ll}
T & 0 \\
0 & 0
\end{array} 
\right]
\q.
$$
Let $Q$ be the partial isometry
$$
Q := 
\left[
\begin{array}{ll}
T & \sqrt{1-TT^*} \\
0 & 0
\end{array} 
\right] 
$$ 
of Proposition \ref{prop1}.
Since $Q$ is a compact perturbation of $\tilde{T}$, 
$$
\ind(Q) =  \ind (\tilde{T}) = \ind (T) 
\q.
$$
Also, since $Q$ is a partial isometry, 
$\dim \kr (Q) = \dim \rang (1-Q^*Q)$ and  
$\dim \coker (Q) = \dim \rang (1-QQ^*)$.  
Hence by Proposition \ref{prop1}, 
\begin{eqnarray*}
\ind (T) = \ind (Q) &:=& 
\dim \kr (Q) - \dim \coker (Q) \\
&=& \dim \rang (1-Q^*Q) - \dim \rang (1-QQ^*)\\
&=& \dim \rang (1-T^*T) - \dim \rang (1-TT^*)\\
&=& - K(T)
\q.
\end{eqnarray*} 
\end{proof}
\noindent
{\bf Remark:}  The above proof that $\ind (T) = 
\dim \rang (1-T^*T) - \dim \rang (1-TT^*)$ is concise and natural
within our context, but may not be the most insightful.
A slightly messier but more straightforward proof 
can be based on the well-known fact 
that for any operator $T$,
the restriction of $ T^*T$ to its initial space 
(defined as the orthogonal complement of its nullspace) 
is unitarily equivalent
to the restriction of $TT^*$ to the closure of {\em its} initial space. 
(The equivalence can be implemented by the 
partial isometry $U$ in the polar decomposition 
$T = U \sqrt{T^*T}$ restricted to its initial space.) 
From this it follows that any nonzero eigenvalue
for $T^*T$ is also an eigenvalue for $TT^*$, 
with the same multiplicity, so that in the expression 
$\dim \rang (1-T^*T) - \dim \rang (1-TT^*)$,
the dimensions of the eigenspaces corresponding to nonzero eigenvalues
cancel, leaving the only contribution to this expression
as $\dim \kr (T) - \dim \kr (T^*) =  \ind (T)$. 
\section{A contraction with non-integral curvature} 
Now we apply Theorem \ref{curvform} to construct a simple example 
of an operator with non-integral curvature, 
in fact with arbitrary real curvature  $\kappa \geq 0$. 
This shows that Theorem \ref{curvthm}'s hypothesis 
that $T$ be pure cannot be omitted.
I have been told that the existence of non-pure contractions 
with non-integral curvatures was implicitly known 
or expected by experts in the field, 
so the interest of the example may lie more in its simplicity
than novelty.  

It is enough to produce a partial isometry $Q$ with  $\rang \, \Delta_Q$ 
spanned by a single unit vector $e$ satisfying 
$$
	\lim_{n \goesto \infty} \| Q^n e \|^2 = \kappa
\q.
$$ 
First suppose $0 \leq \kappa \leq 1$, and set $\lambda := \sqrt{1-\kappa}$. 

Let $T$ be the bilateral weighted shift defined on an orthonormal basis 
$\{e_n\}_{n=-\infty}^\infty$ by:
$$
Te_n := \left\{ 
\begin{array}{ll}
e_{n+1} & \mbox{if $n \neq 0$}\\ 
\lambda e_1 & \mbox{if $n = 0$}
\end{array}
\right.
\q.
$$
Then one routinely computes that $\Delta_T := \sqrt{1-TT^*}$ 
is the rank 1 operator
whose only non-zero eigenvalue is $\sqrt{1-\lambda^2} = \sqrt{\kappa}$, 
with corresponding eigenvector $e_1$.  

Let $Q$ be the associated partial isometry given by Proposition \ref{prop1}.
We may realize $Q$ as acting on a space with orthonormal
basis $\{ e_\infty \} \bigcup \{e_n\}_{n=-\infty}^\infty$
obtained by adjoining a new unit vector named $e_\infty$
to the previous orthonormal basis for the space on which $T$ was defined. 
Then $Q$ is defined by $Qe_n := Te_n$ for $n$ finite, and 
$Qe_\infty := \sqrt{1-\lambda^2} \, e_1 = \sqrt{\kappa} \, e_1 $. 

As in Proposition \ref{prop1}, one routinely computes that 
$\Delta_Q := \sqrt{1-QQ^*}$
is the one-dimensional projection with range spanned by $e_\infty$. 
From Theorem \ref{curvform}
$$
K(Q) = \lim_{n \goesto \infty} \| Q^n e_\infty \|^2 = 
\lim_{n\goesto \infty} \| T^{n-1} \sqrt{\kappa} \,  e_1 \|^2 = \kappa 
\q.
$$ 

This shows that any $\kappa$ with $0 \leq \kappa \leq 1$ can be 
the curvature of some contraction.   
To see that any real number can be the curvature of some contraction,
first check that curvature is additive over direct sums:
for any two contractions $ T_1, T_2$, we have
$$
K(T_1 \oplus T_2) = K(T_1) + K(T_2) 
\q.
$$
This follows routinely from the original definition \re{curvdef} of
curvature, or slightly more easily,
from Theorem \ref{curvform}. 
Then any desired non-negative real curvature can be obtained by direct summing
appropriate copies of the above example. 

\end{document}